\documentclass[reqno]{amsart}

\usepackage{amsmath,amssymb}

\renewcommand{\span}{{\rm Span}}

\newcommand{\eps}{\epsilon}
\newcommand{\e}{{\rm e}}
\newcommand{\sr}{{\rm sr}}

\newcounter{resultcounter}

\newtheorem{lem}[resultcounter]{Lemma}
\newtheorem{prop}[resultcounter]{Proposition}

\newtheorem{rem}[resultcounter]{Remark}

\newcommand{\R}{{\mathbb R}}
\newcommand{\E}{{\mathbb E}}
\renewcommand{\P}{{\mathbb P}}

\begin{document}
 
\title{On the singularity of random matrices with independent entries}

\author{Laurent Bruneau}
\author{Fran\c cois Germinet}
\address{L.B and F.G.: Universit\'e de Cergy-Pontoise, CNRS UMR 8088, Laboratoire AGM, D\'epartement de Math\'ematiques, 95302 Cergy-Pontoise, France.} 
\address{F.G.: Institut universitaire de France.}
\email{laurent.bruneau@u-cergy.fr}
\email{francois.germinet@u-cergy.fr}


\begin{abstract} 
We consider  $n$ by $n$ real matrices whose entries are non-degenerate random variables that are independent but non necessarily identically distributed, and show that the probability that such a matrix is singular is $O(1/\sqrt{n})$. The purpose of this note is to provide a short and elementary proof of this fact using a Bernoulli decomposition of arbitrary non degenerate random variables.
\end{abstract}

\maketitle

\thispagestyle{empty}
\setcounter{page}{1}
\setcounter{section}{1}


\setcounter{section}{0}

\section{Introduction}
Let $M_n=(a_{ij})$ be a random $n\times n$ matrix, where the $a_{ij}$ are independent (non necessarily identically distributed) real random variables. We assume that the r.v. $a_{ij}$ satisfy the following uniform non-degeneracy property
\begin{itemize}
 \item[(H)]   {\em There exists $\rho\in]0,\frac12[$ such that for any $i,j=1,\cdots,n$,  $\P(a_{ij}>x_{ij}^+)>\rho$ and $\P(a_{ij}<x_{ij}^-)>\rho$ for some real numbers $x_{ij}^-<x_{ij}^+$ }.
\end{itemize}

We provide an elementary proof of the following proposition.

\begin{prop}\label{prop} Let $M_n$ be an $n\times n$ matrix whose
  coefficients are independent random variables satisfying \emph{(H)}. Then $\P(M_n$ is singular$)= O(1/\sqrt{n}).$
\end{prop}

The study of the singularity of random matrices goes back, at least, to Koml\'os who showed in \cite{Ko1} that  $\P(M_n$ is singular$)=o(1)$ for independent and identically distributed (iid) Bernoulli entries, namely $a_{ij}=0,1$ with probability $1/2$. Using Sperner's Lemma, Koml\'os noticed that the probability was $O(n^{-1/2})$ \cite{Bo}, a result which has been further extended in \cite{S} to the case of iid entries equally distributed over a finite set. For iid Bernoulli entries, the conjecture is that  $\P(M_n$ is singular$)=(c+o(1))^n$ with $c=\frac12$. Such an exponential behaviour has been successively obtained and improved in \cite{KKS,TV1,TV2} up to $c=\frac34$. The value $c=\frac12$ still seems to be out of reach.

If one turns to general entries, Koml\'os proved in \cite{Ko2} that $\P(M_n$ is singular$)=o(1)$ for independent {\em and identically distributed} non degenerate random variables. Furthermore, as pointed out by Tao and Vu in \cite[Section~8]{TV1}, it follows from their analysis that $\P(M_n$ is singular$)=o(1)$ for independent non degenerate entries, provided Property~(H) holds. Under the same hypothesis Proposition~\ref{prop} asserts that $\P(M_n$ is singular$)=O(n^{-1/2})$.

But the main purpose of this note is to illustrate how the Bernoulli decomposition developed in \cite{AGKW} may be used in order to extend results known for Bernoulli to the general case of independent non degenerate random variables. We perform that illustration by extending Koml\'os's argument as reproduced in \cite{Bo}, to independent random variables satisfying the Property~(H). It is however not clear, at least to the authors, whether results in \cite{TV1,TV2}, and in particular Hal\`asz-type arguments, could be extended in a similar way.

\section{Proof}

Our approach relies on the following lemma which is essentially contained in \cite{AGKW}. For the reader's convenience we sketch its proof in the appendix.

\begin{lem}\label{lemBernoulli}
Let $M_n$ be an $n\times n$ matrix whose
  coefficients are independent random variables satisfying \emph{(H)}. We can decompose the entries of the matrix $M_n$ as follows:
For all $i,j$, there exist two independent random variables $w_{ij}$ and $\eps_{ij}$ and  functions $f_{ij}:]0,1[\to \R$ and $\delta_{ij}:]0,1[\to ]0,+\infty[$ such that 
\\
1. $\eps_{ij}$ is a Bernoulli random variable with parameter $p_{ij}\in ]0,1[$;
\\
2. $w_{ij}$ has the uniform distribution in $]0,1[$; 
\\
3. $a_{ij}=f_{ij}(w_{ij})+\delta_{ij}(w_{ij})\eps_{ij}.$
\\
Moreover, $p_{ij}\in ]1-p_0,p_0[$ for all $i,j$, where $p_0=1-\rho$.
\end{lem}
 
\begin{rem}
It is of crucial importance for us (see (\ref{useSperner}) in the proof of Lemma~\ref{lem:rank}) that $\delta_{ij}>0$. We however do not need here a uniform bound from below on these $\delta_{ij}$. In some situations, one actually does need such a uniform lower bound (see \cite{AGKW}), in which case it is sufficient to modify \emph{(H)} above and require the existence of $x_-<x_+$  independent of $i,j$.
\end{rem}

Thanks to Lemma~\ref{lemBernoulli} and since $w_{ij}$ and $\eps_{ij}$ are independent r.v., we may adopt the following strategy to prove the proposition: 1. do the conditioning with respect to the variables $w_{ij}$, so that, given the $w_{ij}$'s, $M_n$ becomes a sum of a constant matrix and of a random matrix with Bernoulli entries with probabilities $(1-p_{ij},p_{ij})$ and amplitudes $\delta_{ij}(w_{ij})$; 2. estimate, with respect to the Bernoulli variables $\eps_{ij}$, the probability that $M_n$ is
singular following the strategy of \cite{Bo}; 3. take the expectation value with respect to the variables $w_{ij}$.

We shall denote by $\P^{(w)}$ the conditional probability with
respect to the $w_{ij}$ variables, i.e. $\P^{(w)}(\cdot):= \P(\cdot |
\{w_{ij}\}_{i,j} ).$ 

Following \cite{Bo}, we introduce the strong rank of a system of vectors
$S=\{v_1,\cdots,v_n\}$, $\sr(S)$, to be the largest integer $k$ such that any $k$
of the $v_j$'s are linearly independent. For an $m$ by $n$ matrix $A$, we denote by
$\sr_c(A)$ and $\sr_r(A)$ to be, respectively, the strong rank of the system of
columns and of rows of $A$. 

The first ingredient of the proof is the following upper bound on the probability
for an $m$ by $n$ matrix to have a ``not too large'' strong rank.
\begin{lem}\label{lem:strongrank} Let $A$ be an $m$ by $n$ random matrix whose
  coefficients $a_{ij}$ satisfy \emph{(H)}, and $w=(w_{ij})$ be given. Then 
\begin{equation}
 \P^{(w)}(\sr_c(A)<k)\leq
  \left(\begin{array}{c} n \\ k \end{array} \right)
  \frac{p_0^{m-k+1}}{1-p_0} \mbox{ and }\ \P^{(w)}(\sr_r(A)<k)\leq
  \left(\begin{array}{c} m \\ k \end{array} \right)
  \frac{p_0^{n-k+1}}{1-p_0} . 
\end{equation}
\end{lem}

\proof The second statement is clearly equivalent to the first one
(applied to $A^T$). 
By definition of the strong rank, $\sr_c(A)<k$ if and only if there exists $k$ columns of
$A$ which are linearly dependant. It thus suffices to show that for
any $1\leq i_1<\cdots< i_k\leq n$,
\begin{equation}\label{vectrank}
\P^{(w)}(\mathrm{rank}\{v_{i_1},\cdots,v_{i_k}\}<k)\leq \frac{p_0^{m-k+1}}{1-p_0},
\end{equation}
where the $v_j$ denote the columns of $A$. Now we have
\begin{eqnarray}\label{vectrank2}
\lefteqn{\P^{(w)}(\mathrm{rank}\{v_{i_1},\cdots,v_{i_k}\}<k)}\\
 & \leq & \P^{(w)}(v_{i_1}=0)+\sum_{j=1}^{k-1} \P^{(w)}
 (v_{i_{j+1}}\in\span\{v_{i_1},\cdots,v_{i_j}\} | \mathrm{rank}\{v_{i_1},\cdots,v_{i_j}\}=j ).\nonumber
\end{eqnarray}

Let $B$ denote the $m$ by $j$ matrix whose columns are the
vectors $v_{i_1},\cdots, v_{i_j}$. If
$B$ has rank $j$, without loss of generality we may decompose it as 
$B=\left( \begin{array}{c} C \\ D \end{array}\right)$ where $C$ is an
invertible $j$ by $j$ matrix. In a similar way, let us decompose
$v_{i_{j+1}}$ as $v_{i_{j+1}}=\left( \begin{array}{c} Y \\ Z
  \end{array}\right)$, where $Y$ has length $j$ and $Z$
length $m-j$. Note that for $w$ given in the Bernoulli decomposition, the probability of each entry of $Z$ taking a particular value is bounded by $p_0$.

Then, $v_{i_{j+1}}\in\span\{v_{i_1},\cdots,v_{i_j}\}$ iff there
exists a vector $u=(u_1,\cdots,u_j)^T$ such that $Bu=v_{i_{j+1}}$ and
hence iff $Cu=Y$ and $Du=Z$. But since $C$ is
invertible we finally get $v_{i_{j+1}}\in\span\{v_{i_1},\cdots,v_{i_j}\}$ iff $Z=DC^{-1}Y$.
Therefore we have 
\begin{eqnarray}\label{vectspan}
\lefteqn{ \P^{(w)}(v_{i_{j+1}}\in\span\{v_{i_1},\cdots,v_{i_j}\} | \mathrm{rank}\{v_{i_1},\cdots,v_{i_j}\}=j)} \\
& \quad=  \E^{(w)}_Y\left( \P^{(w)}(Z=DC^{-1}Y | Y) \right)
 \leq  \E^{(w)}_Y(p_0^{m-j})=p_0^{m-j},
\end{eqnarray}
where $\E^{(w)}_Y$ denotes the conditional expectation with respect to the variables $w$ and over the random vector $Y$.
Inserting (\ref{vectspan}) into (\ref{vectrank2}) and noting that
$\P^{(w)}(v_{i_1}=0)\leq p_0^m$, this proves (\ref{vectrank}).  
\qed

The second ingredient of the proof is the following improvement of (\ref{vectrank}).
\begin{lem}\label{lem:rank} Let $v_1,\cdots,v_k\in\R^n$ ($k<
  n$) be linearly independent and $X=(a_1,\cdots,a_n)$ a random vector whose
  coefficients satisfy \emph{(H)}. Suppose that $\sr_r(A)=s$ where $A$ is
  the matrix whose columns are the $v_j$'s. Then $\P^{(w)}(X\in \span\{v_1,\cdots,v_k\})\leq Cp_0^{n-k-1}/\sqrt{s}.$
\end{lem}

The above lemma relies on the following generalization of the Littlewood and Offord
problem to the case of non necessarily
identically distributed r.v. and which is an immediate
consequence of an extended version of Sperner's lemma (see
\cite{AGKW}, Lemma 3.2).
\begin{lem}\label{lem:litt-off} 
If $\alpha_1,\cdots, \alpha_s$ are non zero real numbers, $b\in\R$ and $\eps_1,\cdots,\eps_s$ independent Bernoulli random variables with parameter $p_i\in ]1-p_0,p_0[$, then 
$$
\P(\alpha_1\eps_1+\cdots+\alpha_s\eps_s=b)= O(1/\sqrt{s}).
$$
\end{lem}

\noindent{\bf Proof of Lemma~\ref{lem:rank}.}\ \  Let $B$ denote the ($n$ by $k+1$) matrix $A$ augmented with the column vector $X$, and let $r_1,\cdots,r_n$ denote the rows of $B$. If $X\in \span\{v_1,\cdots,v_k\}$ then $B$ has rank $k$, so that without loss of
generality we may assume that $r_1,\cdots,r_k$ are linearly
independent, and that the others $r_j$ depend on these. In particular
$$
\sum_{i=1}^{k+1} \gamma_i a_i=0,
$$
where $\gamma_{k+1}=1$ and, because $\sr_r(B)\geq \sr_s(A)=s$, at least $s$ of the others $\gamma_i$ are
non-zero. Thus, using Lemma \ref{lem:litt-off}, we have, recalling $\delta(w_i)>0$,
\begin{equation}\label{useSperner}
\P^{(w)}\left(\sum_{i=1}^{k+1}\gamma_i a_i=0 \right)
=
\P^{(w)}\left(\sum_{i=1}^{k+1}\gamma_i \delta(w_i) \eps_i= -\sum_{i=1}^{k+1}\gamma_i f(w_i) \right)
\leq C/\sqrt{s+1}.
\end{equation}
Finally, in the same way as in the proof of (\ref{vectspan}), the
$a_i$ for $k+2\leq i\leq n$ are uniquely determined by 
$a_1,\cdots, a_k,$ and thus each of them has a probability at most $p_0$ to take a particular value.
\qed


\noindent {\bf Proof of Proposition~\ref{prop}.} By Lemma~\ref{lemBernoulli} we have 
$$
\P(\mathrm{rank}(M_n)<n)=\E_{\{w_{ij}\}_{i,j}} \left( \P^{(w)}(\mathrm{rank}(M_n)<n \right).
$$
Let $0<\beta<\alpha<1$ to be specified. Let $C_1,\cdots,C_n$ denote the
column vectors of $M_n$ and write $E_k$ for the event that
$C_1,\cdots,C_k$ are linearly independent and $C_{k+1}$ depends on
them. We then have
$$
\P^{(w)}(\mathrm{rank}(M_n)<n)  \leq  \P^{(w)}( \sr_c(M_n)< \alpha
n)+ \sum_{k=\alpha n}^{n-1} \P^{(w)}(E_k).
$$
Indeed, if $\alpha n\leq \sr_c(M_n)<n$ there exists $k\geq \alpha n$
such that $C_1,\cdots,C_k$ are independent but $C_{k+1}$ does depend
on them.

Fix now $\alpha n\leq k<n$, and denote by $A_k$ the $n$ by $k$ matrix whose columns are $C_1,\cdots,C_k$. We
have then 
$$
\P^{(w)}(E_k) \leq \P^{(w)}(\sr_r(A_k)<\beta n)+\P^{(w)}(E_k |\sr_r(A_k)\geq\beta
n). 
$$
Using Lemmas \ref{lem:strongrank} and \ref{lem:rank}, we thus get 
\begin{eqnarray*}
\lefteqn{\P^{(w)}(\mathrm{rank}(M_n)<n)} \\
 & \leq & \left(\begin{array}{c} n \\ \alpha n \end{array} \right)
  \frac{p_0^{n(1-\alpha)+1}}{1-p_0}+\sum_{k=\alpha n}^{n-1} \left( \left(\begin{array}{c} n \\ \beta n \end{array} \right)
  \frac{p_0^{k-\beta n+1}}{1-p_0} +  Cp_0^{n-k-1}/\sqrt{\beta
    n}\right)\\
 & \leq & \frac{p_0}{1-p_0} \left(\begin{array}{c} n \\ \alpha n
   \end{array} \right)p_0^{(1-\alpha) n}+\frac{p_0}{(1-p_0)^2} \left(\begin{array}{c} n \\ \beta n
   \end{array} \right)p_0^{(\alpha-\beta)
   n}+\frac{C}{(1-p_0)\sqrt{\beta n}}\\
 & \leq & C'\left( \e^{n(h(\alpha)+(1-\alpha)\ln p_0)}
   +\e^{n(h(\beta)+(\alpha-\beta)\ln p_0)}+\frac{1}{\sqrt{n}} \right),
\end{eqnarray*}
where $h(x)=-x\ln (x)-(1-x)\ln(1-x)$ is the entropy function and we
made used of the Stirling formula to get the last line. It finally
suffices to take $0<\beta<\alpha<1$ small enough so that
$h(\alpha)+(1-\alpha)\ln p_0$ and $h(\beta)+(\alpha-\beta)\ln p_0$ are
both strictly negative.
\qed

\section{Appendix}

For the reader's convenience, we recall the basic material from \cite{AGKW} and show how to extract from (H) the desired uniform estimates on the Bernoulli decomposition. Namely, we prove Lemma~\ref{lemBernoulli}.

Let $a$ be a random variable satisfying the estimates of Property~(H), with points $x_-<x_+$. We denote by $\mu$ its law and by $F$ its distribution function: $F(u)=\mu(]-\infty,u])$. We set, for any $t\in]0,1[$,
\begin{equation}
G(t):= \inf\{u, F(u)\ge t\}.
\end{equation}
Note that 
\begin{equation}  \label{eq:G2}
G(t)\le u \quad  \Longleftrightarrow \quad F(u)  
\ge t  .
\end{equation}  
We set $p_-=\mu(]-\infty,x_-[)$, $p_+=\mu(]x_+,+\infty[)$, and $p= 1-p_-$. Thanks to (H), $p \ge p_+ > \rho$, and $1-p=p_-> \rho$, so that $p\in]\rho, 1-\rho[$.

Following the idea of \cite[Proof of Theorem~2.1]{AGKW}, define for $t \in ]0,1[$:
\begin{eqnarray} \label{eq:defY}
Y_1(t) &:=& G\left((1-p)t\right) \nonumber \\
Y_2(t) &:=& x_- - G(p_-) + G\left( 1-p + pt\right).
\end{eqnarray}
We always have $G(p_-)\le x_-$, but these two numbers may differ. Note that since $\mu(]-\infty,x_-[)=p_-$, if it turns out that $G(p_-)<x_-$, then for any $u\in]G(p_-),x_-[$ one has $ p_-\le F(u) \le p_-$, and thus $\mu(]G(p_-),x_-[)=0$ (this is the reason for the definition of $Y_2$ which differs from the one of \cite[Proof of Theorem~2.1]{AGKW}). We then let
\begin{eqnarray} 
     f(t)  & := &  Y_1(t),    \\
      \delta(t)  & := &  Y_2(t) - Y_1(t)  \label{eq:delta} ,
\end{eqnarray}
so that, if $\eps$ is a Bernoulli variable with probabilities $(1-p,p)$ and $t$ a random variable with uniform distribution in $]0,1[$, we do have
\begin{equation} 
a = f(t) + \delta(t)\eps .
\end{equation}  
It remains to prove that $\delta(t)>0$ almost surely.
Obviously, combining (\ref{eq:defY}) and (\ref{eq:delta}), $\delta(t)\ge x_- - G(p_-)$ for any $t$, so that if $x_->G(p_-)$ we are actually done.

Suppose $G(p_-)=x_-$. We claim that $T_1=1>T_2$ where
\begin{eqnarray}
 T_1   =  & \inf \{ t\in ]0,1[ \, : \, Y_1(t) = G(p_-) \}  \qquad &
\mbox{(arrival time of $Y_1$)}  \, , \nonumber \\
 T_2   = & \sup \{ t\in ]0,1[ \, : \, Y_2(t) = G(1-p + 0 ) \}  \qquad &
\mbox{(departure time of $Y_2$)} \,  .   \nonumber 
\end{eqnarray}
As in \cite[Proof of Theorem~2.1]{AGKW}) the latter then implies that $\delta(t)>0$ for all $t$. It is easy to see that $T_2\le (1-p_- -p_+)/(1-p_-) < 1$. It remains to show that $T_1=1$. Suppose $T_1<1$. For any $t\in]T_1,1[$ and for any $u<x_-$, one has $x_-=G(p_- t)>u$. Then (\ref{eq:G2}) implies that $F(u) < p_- t$, and thus we get the following contradiction
\begin{equation} 
p_- = \mu(]-\infty,x_-[) = \sup_{u<x_-} F(u) \le p_-t < p_-.
\end{equation}

\section*{Acknowledgement}
It is the pleasure of the authors to thank Abel Klein for his comments on a preliminary version of this note. F.G. also thanks Michael Aizenman and Simone Warzel for their warm hospitality in Princeton.

\end{document}